\documentclass[reqno]{amsart}

\newcommand{\RR}{\mathbb{R}}
\newcommand{\ZZ}{\mathbb{Z}}

\newcommand{\eb}{\mathbf{e}}
\newcommand{\nb}{\mathbf{n}}
\newcommand{\qb}{\mathbf{q}}
\newcommand{\pb}{\mathbf{p}}
\newcommand{\ub}{\mathbf{u}}
\newcommand{\xib}{\boldsymbol{\xi}}
\newcommand{\etab}{\boldsymbol{\eta}}
\newcommand{\zetab}{\boldsymbol{\zeta}}
\usepackage{flexisym}
\usepackage{amssymb}
\usepackage{graphicx}
\usepackage{xcolor}
\usepackage{bm}
\usepackage{amsmath}
\usepackage{tikz}
\usepackage{algorithm}
\usepackage[noend]{algpseudocode}
\usepackage{caption}
\usepackage{subcaption}
\newcommand*\Let[2]{\State #1 $\gets$ #2}
\usepackage{ amssymb }
\usepackage{geometry}
\usepackage{tabularx}

\title{Strong Convergence of Integrators for Nonequilibrium Langevin Dynamics}
\date{\today}
\author{Matthew Dobson and Abdel Kader Geraldo}
\address{Department of Mathematics and Statistics, 710 N Pleasant St., Amherst, MA 01003}
\begin{document}
\begin{abstract}
Several numerical schemes are proposed for the solution of
Nonequilibrium Langevin Dynamics (NELD), and the rate of convergence is
analyzed.    Due to the special deforming boundary conditions used, care must
be taken when using standard stochastic integration schemes, and we demonstrate
a loss of convergence for a naive implementation.  We then present several
first and second order schemes, in the sense of strong convergence.
\end{abstract}

\maketitle

\section{Introduction}
Nonequilibrium molecular dynamics techniques are employed in the study of
microscopic systems undergoing steady, nonconstant flow, for example, in
the study of polymer melts.  A wide range of dynamics,
including both deterministic and stochastic equations have been 
proposed for such simulations~\cite{Evans, Snook}.  

We examine the rates of strong convergence of several numerical methods for the
simulation of Nonequilibrium Langevin Dynamics (NELD)~\cite{McPhie, Snook,
dlls}.  Let $\qb,\pb \in \RR^{3 N}$ denote the positions and 
velocities of a set of particles, then NELD is given by
\begin{equation} \label{neld}
\begin{split}
d \qb &= \pb \,  d t\\
d \pb &=    ( - \Call{$\nabla E$}{\qb}  -\gamma (\pb-A \qb) + A \pb) \, dt + \sigma \, dW
\end{split}
\end{equation}
where  $- \nabla E(\qb)$ are the interparticle forces, $ W$ is a standard
$3N$-dimensional Brownian motion, $ \sigma $ and $ \gamma $ are scalar
constants satisfying the fluctuation-dissipation relation 
\begin{equation}
\label{eq:fluc_diss}
\begin{split}
\gamma  &=  \frac{1}{2} \sigma^2 \beta
\end{split}
\end{equation}
where $\beta= \frac{1}{k_B T}$ is the inverse temperature, 
and $A \in \RR^{3N \times 3N}$ is trace-free block diagonal linear background flow
matrix.  The diagonal entries of $A$ are identical $3 \times 3$ trace-free 
diagonal matrices, corresponding to the macroscale background flow $A = \nabla \ub.$   

To simulate the bulk motion of particles with a mean background flow $A,$
specialized periodic boundary conditions are employed, in particular, a
particle with the coordinates $(\qb,\pb)$ has periodic images at $(\qb+L_t
\nb,\pb+A L_t\nb)$, where $L_t : [0,\infty) \rightarrow \RR^{3N \times 3N}$
is a block diagonal matrix whose $3 \times 3$ identical blocks 
denote the matrix of lattice basis vectors at time $t$ and $\nb \in
\mathbb{Z}^{3N}.$  The images of a single particle do not have the same velocity,
rather they are consistent with the mean flow, and this in turn implies that
the periodic lattice generated by $L_t$ deforms with the flow.
Care is needed to ensure that the lattice does not become degenerately deformed
where the minimum replica distance goes to zero.  Techniques have been
developed in the papers~\cite{Lees-Edwards,Kraynik,Dobson,Hunt} which choose
initial lattice vectors $L_0$ such that the minimum replica distance in the lattice stays bounded away
from zero, and the simulation box is remapped so that the geometry stays
regular.   We will consider the Generalized Kraynik-Reinelt (GenKR) boundary
conditions developed in~\cite{Dobson,Hunt}, which can handle general
three-dimensional incompressible flows.

In this paper, we will focus on the strong convergence properties of certain
common stochastic integrators applied to NELD, seeing how the periodic boundary
conditions interact with the convergence.  In particular, we will see that a
naive implementation of certain standard schemes show a breakdown in
convergence due to the interaction of the integrator with the GenKR boundary
conditions.  We will then develop schemes that avoid this convergence problems
and compute the order of convergence by using the Ito-Taylor expansion.
Several standard first and second order schemes will demonstrated numerically
and analytically.

\section{Ito-Taylor  expansion of the nonequilibrium Langevin dynamics}

In this section, we compute the Ito-Taylor expansion for the NELD up to 
second order, which will be used in the error analysis of the numerical 
integrators.  We also set the notation for the application of boundary
conditions as the motion of the replicas plays an important role in the 
analysis of the numerical schemes.  

\subsection{Review of the Ito-Taylor  expansion }

We express the NELD~\eqref{neld} in integral form, 
\begin{equation} \label{sde_integral}
\begin{split}
X(t)=X(t_0)+ \int_{t_0}^{t}  \Call{C}{X(s)} \, ds +\int_{t_0}^{t} \Sigma \,
dW(s)
\end{split}
\end{equation}
where 
$$C\left(\begin{bmatrix} \qb \\ \pb \end{bmatrix} \right)=
\begin{bmatrix}
\pb  \\
- \nabla E(\qb)  -\gamma (\pb-A \qb) + A \pb  
\end{bmatrix} \qquad
\Sigma=\begin{bmatrix}
0 & 0\\
0 & \sigma I
\end{bmatrix}.$$
The Ito formula for a scalar-valued function $G(X(t))$ of the solution $X(t)$ is given by 
\begin{equation} \label{eq1}
\begin{split}
G(X(t))=G(X(t_0))+ \int_{t_0}^{t} L^0 G(X(s)) \, ds 
                     + \int_{t_0}^{t} L^1 G(X(s)) \, dW(s)
\end{split}
\end{equation}
where the operators $L^0$ and $L^1$ are given by:
\begin{equation} \label{eq2}
\begin{split}
L^0&=\Call{C}{} \cdot \nabla_x 
+\frac{1}{2}\Sigma \Sigma^{T} : \nabla_x^2,\quad
L^1=\Sigma \nabla_x \, . 
\end{split}
\end{equation}
Over a small time interval $\Delta t := t - t_0,$ we apply the Ito formula to 
equation~\eqref{sde_integral} and expand up to second order, noting that several terms
cancel due to the form of $C$ and $\Sigma$, arriving at 
\begin{equation} \label{main1}
\begin{split}
X(t)&=X(t_0)+  \Call{C}{X(t_0)}\int_{t_0}^{t} ds+\Sigma \int_{t_0}^{t}dW(s) + L^0  \Call{C}{X(t_0)} \int_{t_0}^{t}\int_{t_0}^{s}  du \, ds\\
&\quad+ L^1  \Call{C}{X(t_0)} \int_{t_0}^{t}\int_{t_0}^{s}  dW(u) \, ds
+R 
\end{split}
\end{equation}
where the remainder of order O($\Delta t^{5/2 }$) term is given by 
\begin{equation} \label{eq3}
\begin{split}
R&=\int_{t_0}^{t}\int_{t_0}^{s}\int_{t_0}^{u} L^1 L^0 C(X(v))  \, dW(v) \, du \, ds
+ \int_{t_0}^{t}\int_{t_0}^{s}\int_{t_0}^{u} L^0 L^0  C(X(v))\, dv \, du \, ds.
\end{split}
\end{equation}
We recall the following facts about the covariance of $W(t)$ and its integral,
which are useful in developing numerical schemes~\cite{kloeden}: 
\begin{equation} \label{eq4}
\begin{split}
\mathbb{E}(W_i(s)W_j(s^')) &=\delta_{ij}\min(s,s^') \\
\mathbb{E}[(W_i(s)-W_i(s^'))(W_j(s)-W_j(s^'))]&=\delta_{ij} \Delta t\\
\mathbb{E}\left[(W_i(t+\Delta t)-W_i(t))\int_{t}^{t+\Delta t}(W_j(s)-W_j(t)) \, ds\right]&=\frac{1}{2}\delta_{ij} \Delta t^2\\
\mathbb{E}\left[\int_{t}^{t+\Delta t}(W_i(s)-W_i(t)) \, ds \int_{t}^{t+\Delta t}(W_j(s)-W_j(t)) \, ds \right]&=\frac{1}{3}\delta_{ij} \Delta t^3
\end{split}
\end{equation}
Therefore, truncating the expansion of NELD to second order and letting $(\qb, \pb)$ denote the 
coordinates at time $t_0,$ we arrive at
\begin{gather}\label{eq:ito-tay}
\begin{split}
\begin{bmatrix} \qb(t) \\ \pb(t) \end{bmatrix}
&=
\begin{bmatrix} \qb \\ \pb\end{bmatrix}+
\begin{bmatrix}
\pb  \\
\Call{F}{\pb,\qb}  
\end{bmatrix}\Delta t
+  \Delta t^{1/2} 
\begin{bmatrix}
0&0\\
0&\sigma I
\end{bmatrix} \etab \\
&\quad +
\begin{bmatrix}
\Call{F}{\pb,\qb}\\
(- \Call{$\nabla^2 E$}{\qb} +\gamma A) \pb +(A-\gamma I)\Call{F}{\pb,\qb} 
\end{bmatrix} \frac{\Delta t^2}{2} \\
&\quad +\sigma \Delta t^{3/2} \begin{bmatrix}
0 & I\\
0 &  (A-\gamma I)\\
\end{bmatrix} \left(\frac{1}{2}\etab+\frac{1}{2\sqrt{3}}\zetab\right) 
+ \Call{O}{ \Delta t^{5/2}}
\end{split}
\end{gather}
where $  \Call{F}{\pb,\qb} =- \Call{$\nabla E$}{\qb}  -\gamma (\pb-A\qb) + A \pb $ 
and where 
\begin{equation*}
\begin{split}
\etab &= \frac{W(t) - W(t_0)}{\Delta t},  \\
\zetab &= \frac{2 \sqrt{3}}{\Delta t^{3/2}} \int_{t_0}^t \left( W(s) - W(t_0) \right) \, ds \,  - \, \sqrt{3} \, \etab.
\end{split}
\end{equation*}  Note the scaling of 
stochastic terms has been chosen so that $\etab, \zetab \sim \mathcal{N}(0,1)$ are independent Gaussian random variables.

\subsection{Nonequilibrium Boundary Conditions}
When analyzing the truncation error for the scheme, it is important to account
for the nonequilibrium periodic boundary conditions, particularly the fact that
replicas do not all have the same velocity.  A particle with coordinates $(\qb,
\pb)$ has periodic images at $(\qb + L_t \nb, \pb + A L_t \nb).$  The NELD
equations are invariant under a translation of the system by choosing new $\nb.$
In particular, $d \qb = \pb \, dt$ holds for all particle images, so that 
$$\frac{d}{dt}(\qb+L_t \nb)=\pb+A L_t \nb $$ which  imply that  the simulation
box deforms with the flow, $ \frac{d}{dt} L_t=AL_t$ with a solution $L_t=\exp (A
t)L_0$. 

During a simulation step, one or more particles can leave the simulation box,
whereupon they are remapped in accordance to the periodic boundary conditions.
This can also be viewed as no longer tracking the position of the particles
that started at $(\qb, \pb),$ but tracking the particles at $(\qb + L_t \nb, \pb + A L_t \nb)$
for some $\nb \in \ZZ^{3N}.$  We show in the following that the timing
of applying the periodic boundary conditions affects the rate of convergence
for the numerical scheme, in fact, reducing the strong rate of convergence 
below first order for a pair of schemes.  When computing the local truncation
error, we compare the final position after the numerical step and periodic 
remapping with the Taylor-Ito expansion of the corresponding replica, which
may have started outside the simulation box.  That is, if we are now tracking
the particle at $(\qb(t) + L_t \nb, \pb(t) + A L_t \nb),$ we compare with the
particle that started at $(\qb(t_0) + L_{t_0} \nb, \pb(t_0) + A L_{t_0} \nb).$ 
We note that the
Taylor-Ito expansion will now have terms from the deformed lattice
vectors,
$$ \qb(t) + L_t \nb = \qb(t) + L_{t_0} \nb + \Delta t A L_{t_0} \nb + \frac{1}{2} \Delta t^2 A^2 L_{t_0} \nb
+ O(\Delta t^3),
$$ 
with a similar expression for $\pb(t) + A L_t \nb.$

\section{Strong Convergence and a Numerical Experiments}  
For a stochastic process, there are several notions of convergence one 
can consider, including strong convergence, weak convergence, or 
convergence of the dynamics to an invariant measure.  In the following,
we consider strong convergence of the proposed numerical schemes.  Given
a stochastic process $X(t),$ we say that the numerical method generating $X_h(t)$ 
has strong order of convergence $r$ if
$$
\mathbb{E}( |X(t) - X_h(t) |) \leq C h^r
$$
for some $C>0$ and all sufficiently small $h>0.$

For each of the described algorithms, we perform a benchmark test to numerically
compute the strong rate of convergence.  For each time, we numerically estimate
the rate of convergence in the $\ell^2$ norm of both the position $\qb$ and 
momentum $\pb,$ and in each case the convergence is observed to behave similarly
in $\qb$ and $\pb.$  We simulate a system of $1728$ particles, having 
the Weeks-Chandler-Anderson (WCA) interparticle interaction potential,
\begin{equation*}
\phi(r) = \begin{cases} \frac{1}{r^{12}} - \frac{1}{r^6} + \frac{1}{4} & r < 2^{1/6} \\
0 & r \geq 2^{1/6},
\end{cases}
\end{equation*}
which is purely repulsive and continuously differentiable.
The background flow used for all tests is a uniaxial extensional flow,
whose diagonal blocks are given by
\begin{equation*}
A = \begin{bmatrix}
0.2 & 0 & 0 \\
0 & -0.1 & 0 \\
0 & 0 & -0.1
\end{bmatrix}
\end{equation*}

In Table~\ref{tab:param}, we list the parameters used for the numerical experiments.  
\begin{table}
\begin{tabular}{ |p{5cm}|p{1.5cm}||p{5cm}|p{1.5cm}| }
\hline
Parameter & Value &  Parameter & Value\\
\hline
Time step                 ($\Delta t$) & 0.000025 &  Friction coefficient    $ \gamma$ & 1.0   \\
Simulation time (T)            & 1.0 &   Inverse temperature    $ \beta$ & 1.0     \\
Number of Particles    & 1728  & Simulation Box Side Length & 15    \\
\hline
\end{tabular}
\caption{\label{tab:param}List of parameters used for the convergence tests of the 
nonequilibrium algorithms}
\end{table}
To gather statistics, we average 200 runs for each numerical experiment.  Each run
is initialized by first running an equilibrium simulation using standard Langevin dynamics, 
which acts to draw
the initial condition according to the Gibbs measure corresponding to equilibrium.
Then, at time zero, the background flow is turned on, so that the system evolves
from the initial state according to NELD equations of motion, including deformation
of the simulation box.  For a given initial state, the nonequilibrium simulation
is run with five different stepsizes: $\Delta t, 2 \Delta t, 4 \Delta t, 8 \Delta t,$ and $16 \Delta t.$
To measure the strong convergence, the Brownian motion for each simulation is the same, 
see~\cite{higham} for an introduction
to numerical computation of stochastic order of convergence.  We then compute and
report differences in the $\ell^2$ norm of the system, $\eb_h(t) = \| \qb_h(t) - \qb_{2h}(t) \|_2$
and estimate the order of convergence by
$$
\mathrm{ord}(t) = \frac{ \log( \eb_{2h}(t) ) - \log(\eb_{h}(t) )}{ \log 2}.
$$

\section{Failed Schemes}
 We recall that the Euler-Maruyama
scheme has order one-half when applied to stochastic equations with multiplicative
noise, but order one equations with additive noise, and we will show
in the next section that it also has order one for the NELD case. In this section, we will analyze two schemes which converge to first order
when applied to  equilibrium Langevin dynamics ($A = 0$), but which fail
to converge to that order in the NELD case. Both schemes will be modified to
have first order in the following section.

\subsection{Symplectic Euler A (SE-A)}
For deterministic Hamiltonian dynamics, there are two types of 
Symplectic Euler integrators, either the position is integrated first
then the momentum in Symplectic Euler A (SE-A) or momentum then position in 
Symplectic Euler B (SE-B).  For the NELD case, we modify the SE-A scheme by applying  
the PBCs at the beginning of the computation and  after  incrementing the 
position, then integrating the stochastic terms using an explicit step,
as in in the Euler-Maruyama scheme.  The SE-A algorithm is described as follows:
\begin{algorithm}[H]\small
\caption{Symplectic Euler A (SE-A)
\label{SEA1}}
\begin{algorithmic}
\Statex
\For{$k = 1 \dots$ Nsteps}
\State \Call{GenKR}{$\qb^k,\pb^k$}  \Comment{Apply PBCs}
\Let{$ \qb^{k+1} $}{$  \qb^k + \pb^k \Delta t$}
\State \Call{GenKR}{$\pb^k,\qb^{k+1}$}                \Comment{Apply PBCs}
\Let{$ \pb^{k+1} $}{$  {\pb}^k +$  \Call{F}{$ {\pb}^k, {\qb}^{k+1}$}$ \Delta t+\sigma \sqrt{\Delta t} \etab$}
\EndFor
\end{algorithmic}  
\end{algorithm}
The particles are potentially remapped twice during the algorithm.  We translate
the pseudocode to our numerical scheme, giving
\begin{equation*}
\begin{split}
&\qb_1 = \qb^k - L \nb_1, \qquad  \pb_1 = \pb^k - A L \nb_1 \\
&\qb_2 = \qb_1 + \pb_1 \Delta t \\
&\qb^{k+1} = \qb_2 - L \tilde{\nb}_2, \qquad  \pb_2 = \pb_1 - A L \tilde{\nb}_2 \\
&\pb^{k+1} = \pb_2 + F(\pb_2, \qb^{k+1}) \Delta t + \sigma \sqrt{\Delta t} \, \etab 
\end{split}
\end{equation*}
where $F(\pb, \qb) := - \nabla E(\qb) - \gamma (\pb - A \qb) + A \pb.$
Letting $\nb_2 = \nb_1 + \tilde{\nb}_2$ and expanding the algorithm, we get
\begin{equation}
\begin{split}
\begin{bmatrix} \qb^{k+1} \\ \pb^{k+1} \end{bmatrix}
&=
\begin{bmatrix} \qb^{k}-L \nb_2 \\ \pb^{k}-AL\nb_2 \end{bmatrix}+
\begin{bmatrix}
\pb^k -AL\nb_1 \\
\Call{F}{\pb^k -AL\nb_2,\qb^k -L\nb_2+(\pb^k -AL\nb_1 ) \Delta t}  
\end{bmatrix} \Delta t+
\begin{bmatrix}
0&0\\
0&\sigma I
\end{bmatrix} \etab\\
&=
\begin{bmatrix} \qb^{k}-L\nb_2 \\ \pb^{k}-AL\nb_2 \end{bmatrix} 
+
\begin{bmatrix}
\pb^k -AL\nb_1 \\
\Call{F}{\qb^k,\pb^k} -A^2L \nb_2  
\end{bmatrix} \Delta t
+
\begin{bmatrix}
0&0\\
0&\sigma I
\end{bmatrix} \sqrt{\Delta t} \etab\\
&
\quad +
\begin{bmatrix}
0 \\
\gamma A(\pb^k-AL\nb_1) - \Call{$\nabla^2 E$}{\qb^k}(\pb^k -AL \nb_1 )
\end{bmatrix} \Delta t^2.
\end{split}
\label{eq:sea_exp}
\end{equation}
Comparing this to~\eqref{eq:ito-tay} applied to the replica that began at $(\qb(t_0) + L_{t_0} \nb_2, \pb(t_0) + A L_{t_0} \nb_2),$ we
find that the leading order terms in the truncation error are
\begin{equation*}\label{q1}
\begin{split}
T_1 &=
\begin{bmatrix}  A L \tilde{\nb}_2\\ 0\end{bmatrix}\Delta t+ \begin{bmatrix}
\Call{F}{\pb,\qb} + A^2 L \nb\\
(\Call{$\nabla^2 E$}{\qb} -\gamma A) (\pb - 2 A L \nb_1) +(A-\gamma I)\Call{F}{\pb,\qb} 
\end{bmatrix} \frac{\Delta t^2}{2} \\
&\quad +\sigma \Delta t^{3/2} \begin{bmatrix}
0 & I\\
0 &  (A-\gamma I)\\
\end{bmatrix} \left(\frac{1}{2}\etab+\frac{1}{2\sqrt{3}}\zetab\right) 
\end{split}
\end{equation*}
Note that since the local truncation error contains a term of $O(\Delta t),$ we
do not expect to have first-order convergence.  However, this term is only non-zero when
there is particle motion across the boundary, so convergence is still possible.  We see in the following numerical
experiment that for the chosen parameters, the order is reduced but the scheme is convergent.
This term arises due to the application of periodic boundary conditions between
the $\qb$ update and $\pb$ update.
 
\subsubsection*{Numerical Result}
\begin{figure}[h]
\begin{subfigure}{.5\textwidth}
\centering
\includegraphics[width=70mm]{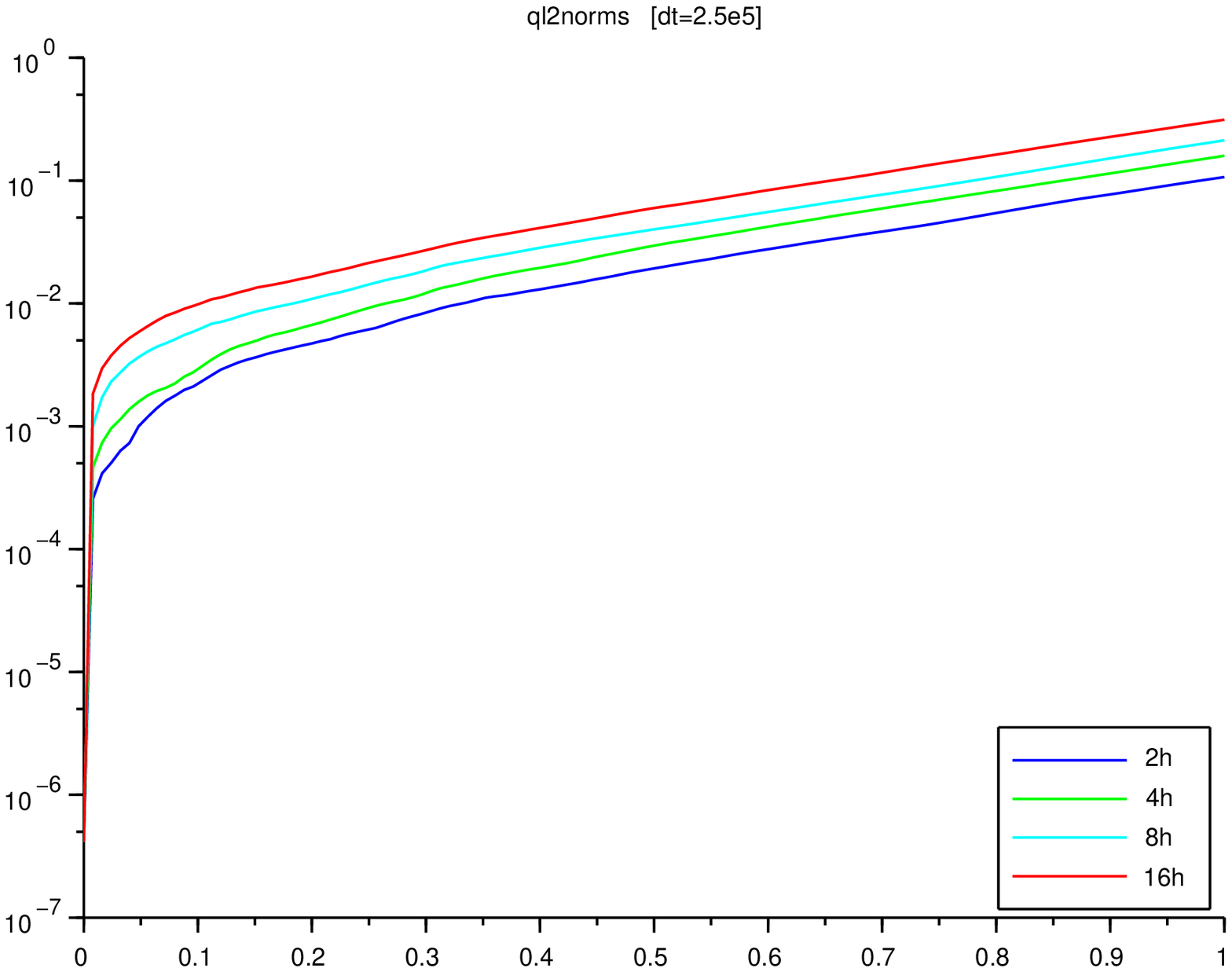}
\end{subfigure}%
\begin{subfigure}{.5\textwidth}
\centering
\includegraphics[width=70mm]{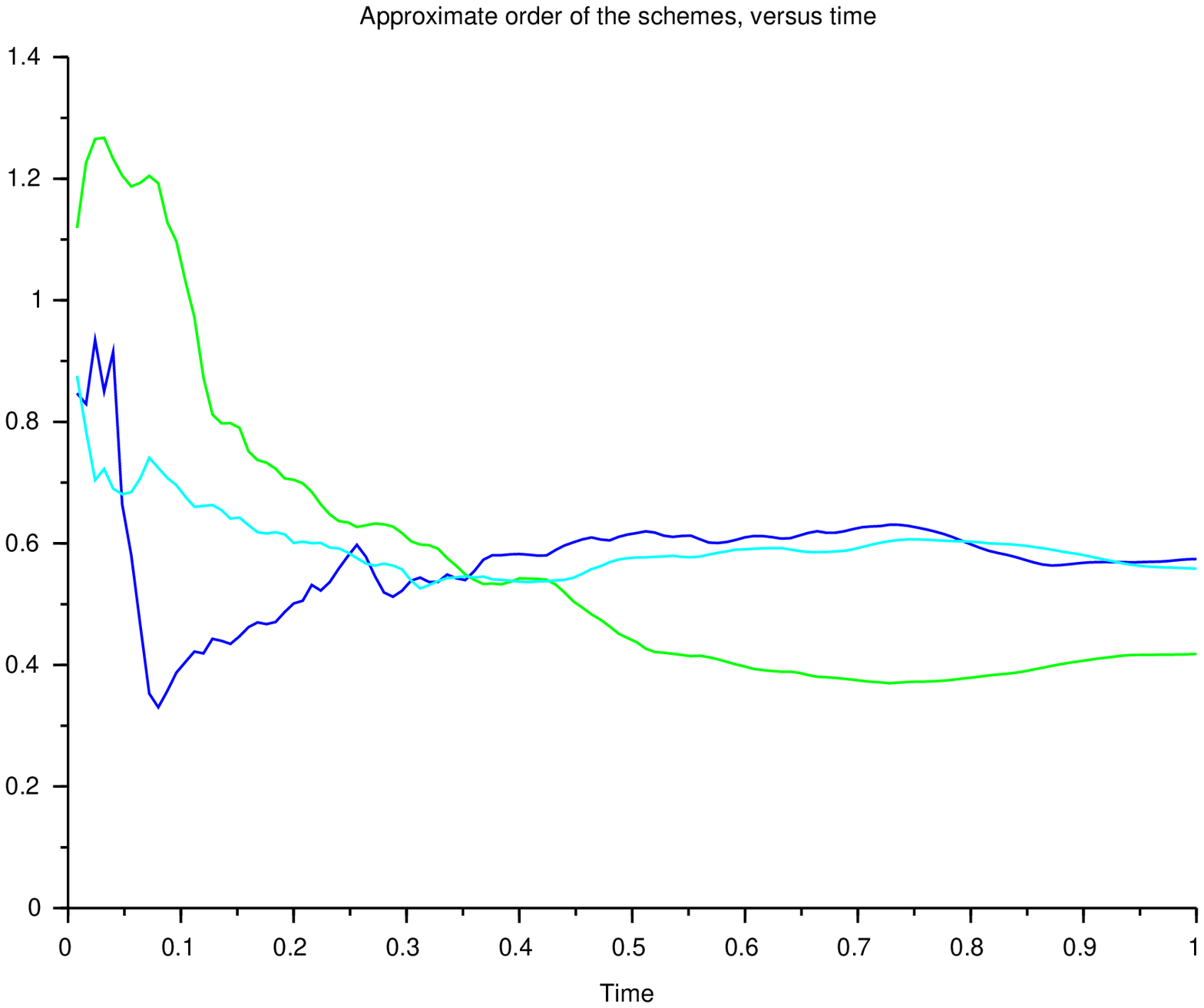}
\end{subfigure}%
\caption{\label{fig:sea}Figure showing the error and order of convergence with the  SE-A scheme}
\end{figure}

Figure~\ref{fig:sea} plots both the errors $\eb_h(t)$ as well as the
numerically estimated order $\textrm{ord}_h(t)$ for the SE-A implementation.
The plots show the lack of first-order convergence for the chosen parameters,
and it can be seen in the graph that the scheme converges approximately at
order 1/2, though this is both irregular among the various runs and depends on
the chosen parameters.  We correct this problem in Section~\ref{sec:ordone} by
delaying the application of PBCs in the algorithm, arriving at the SE-AC
(Symplectic Euler-A, corrected) scheme. 

\subsection{ABAPO}
We consider now a splitting scheme, where the Ornstein-Uhlenbeck portion is
integrated analytically.  We use the terminology from the splitting-scheme
framework of~\cite{Leimkuhler}, where we split the NELD dynamics into three
portions,
\begin{equation}\label{BAO1}
\begin{split}
d\begin{bmatrix} \qb \\ \pb \end{bmatrix}
&=
\underbrace{\begin{bmatrix} \pb \\ 0\end{bmatrix} \, dt}_{A}+
\underbrace{\begin{bmatrix} 0 \\ - \Call{$\nabla E$}{\qb}\end{bmatrix} \, dt}_{B}+
\underbrace{\begin{bmatrix} 0 \\  A \pb\end{bmatrix} \, dt}_{P}+
\underbrace{\begin{bmatrix} 0 \\ -\gamma (\pb-A \qb)  \end{bmatrix}  \, dt+ \begin{bmatrix} 0 \\ \sigma \, dW \end{bmatrix} }_{O}
\end{split}
\end{equation}
each of these split portions can be analytically integrated, $A,$ $B,$ and $P$ trivially so,
while the exact solution of the $O$ part is
\begin{equation}\label{fff2}
\begin{split}
\pb(t)= \exp(-\gamma \Delta t) \pb(t_0)+(1-\exp(-\gamma \Delta t)) A \qb(t_0)+ \sqrt{ \beta^{-1} (1-\exp(-2 \gamma \Delta t)} \xib,
\end{split}
\end{equation}
where 
\begin{equation*}
\xib =  \left(\frac{\beta}{1-\exp(-2 \gamma \Delta t }\right)^{1/2}  \exp(-\gamma t) \int_{t_0}^t \sigma \exp(\gamma s) \, dW(s) \sim \mathcal{N}(0,1).
\end{equation*}
We choose an ABAPO splitting, which is used in~\cite{Burrage}, to arrive
at the numerical integrator 
\begin{equation}\label{ABAPO22}
\begin{split}
\psi_{ABAPO}^{\Delta t}=\exp\left(\frac{\Delta t}{2}\mathcal{L}_A\right)
\exp(\Delta t\mathcal{L}_B)\exp\left(\frac{\Delta t}{2}\mathcal{L}_A\right)
\exp(\Delta t \mathcal{L}_O)
\end{split}
\end{equation}
where $ \mathcal{L}_f $ is the corresponding operator for the vector field $f$.  
It may be noticed that the phase ABA is the
standard Verlet method. 

\begin{figure}[ht]
\begin{subfigure}{.5\textwidth}
\centering
\includegraphics[width=70mm]{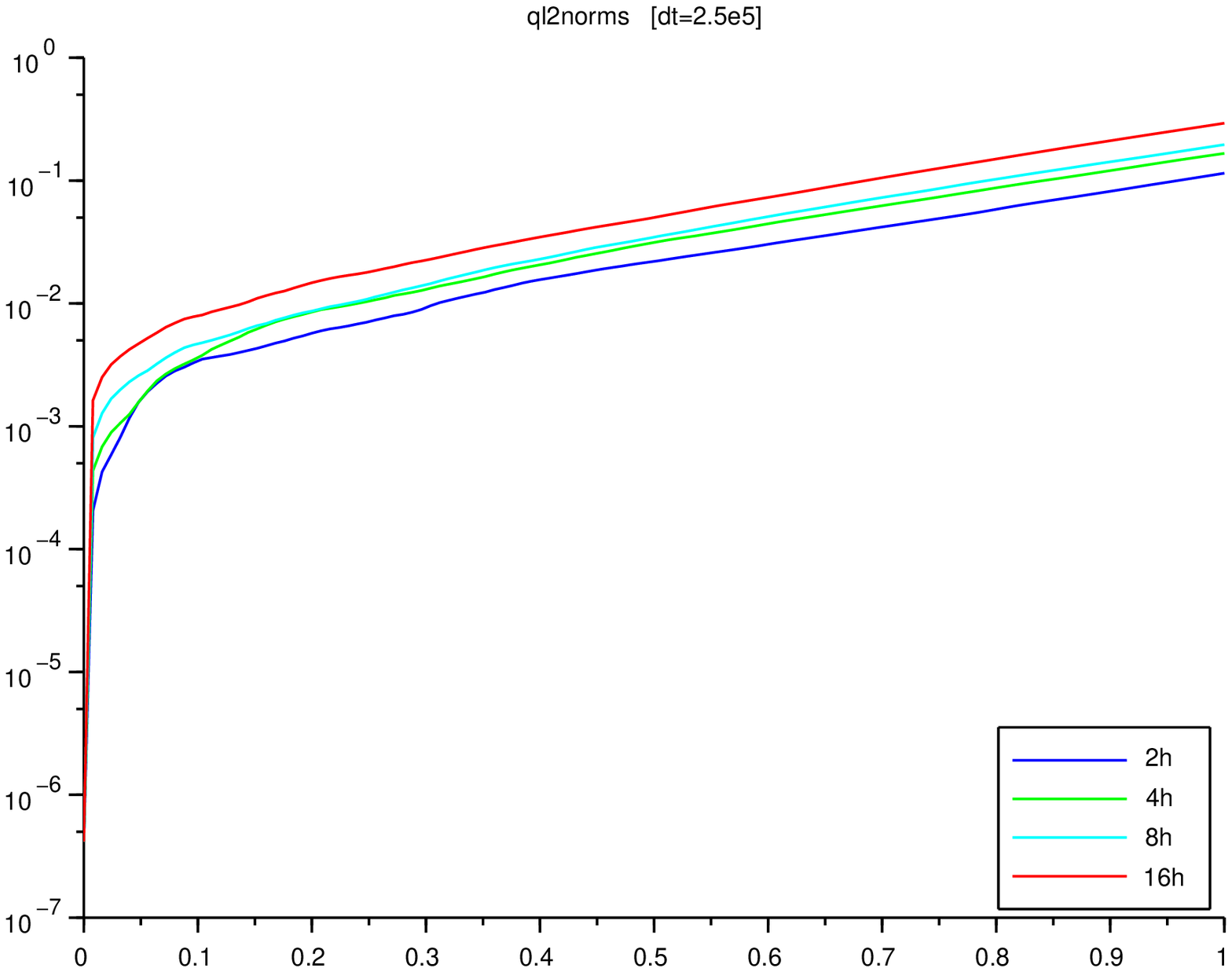}
\end{subfigure}%
\begin{subfigure}{.5\textwidth}
\centering
\includegraphics[width=70mm]{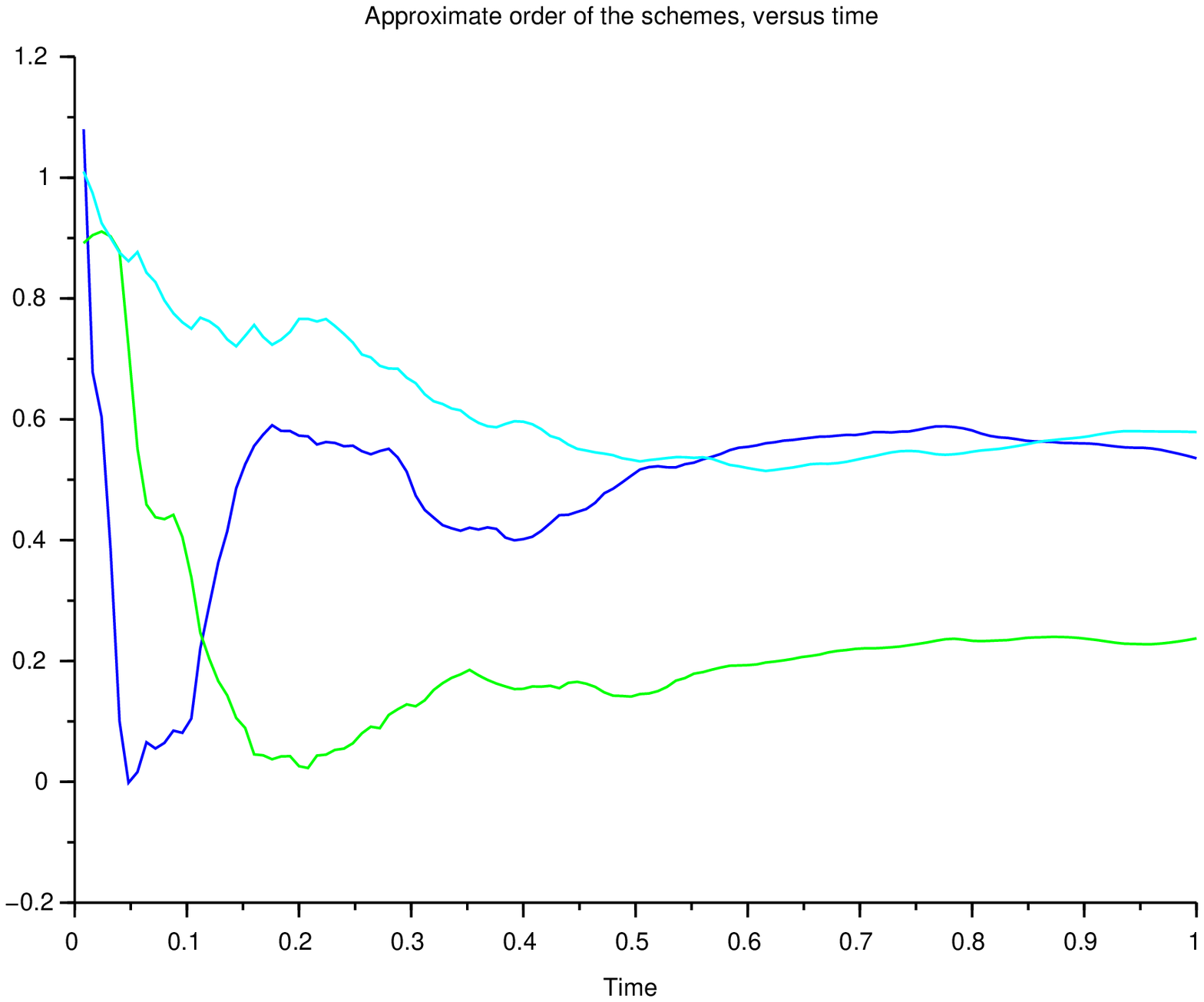}
\end{subfigure}%
\caption{
\label{fig:ABAPO}Figure showing the error and order of convergence with the  ABAPO scheme}
\end{figure}

 GenKr PBCs is applied at the
beginning of the scheme and after each integration in the position, in order to
keep the particle inside the box.
The ABAPO algorithm is described as follows:
\begin{algorithm}
\caption{ABAPO
\label{ABAPO}}
\begin{algorithmic}
\For{$k \gets 1 \dots $ Nsteps}
\State \Call{GenKR}{$\qb^k,\pb^k$} \Comment{Apply PBCs}
\Let{$ \pb^{k+\frac{1}{4}} $}{$  \pb^k + \frac{\Delta t}{2} \Call{$\nabla E$}{\qb^k}  $}
\Let{$ \qb^{k+1} $}{$  \qb^k +  \frac{\Delta t}{2} \pb^{k+\frac{1}{4}} $}
\State \Call{GenKR}{$ \pb^{k+\frac{1}{4}},\qb^{k+1}$}               
\Let{$  \pb^{k+ \frac{1}{2}} $}{$\pb^{k+ \frac{1}{4}} + \frac{\Delta t}{2} \Call{$\nabla E$}{ \qb^{k+1}} $}
\Let{$ \tilde{ \pb^{k+ \frac{1}{2}}}   $}{$ \exp(\Delta t A) \pb^{k+\frac{1}{2}} $}
\Let{$ \pb^{k+1} $}{$ \gamma \tilde{ \pb^{k+ \frac{1}{2}}}  +(1-\gamma)A \qb^{k+1} +\sqrt{ \beta^{-1} (1-\exp(-2 \gamma \Delta t)} \xib$}
\EndFor
\end{algorithmic}
\end{algorithm}
We observe a similar loss of order of convergence as in the SE-A case, as displayed 
in Figure~\ref{fig:ABAPO}.  There doesn't seem to be a single observed rate of
convergence, though it is clear that it is lower than first order.

\section{First Order NELD Algorithm}
\label{sec:ordone}
In this section, we will analyze four first order NELD schemes, two of which
are corrected versions of Algorithms~\ref{SEA1} and~\ref{ABAPO}. 

\subsection{Euler-Maruyama}

\begin{figure}[h]
\begin{subfigure}{.5\textwidth}
\centering
\includegraphics[width=70mm]{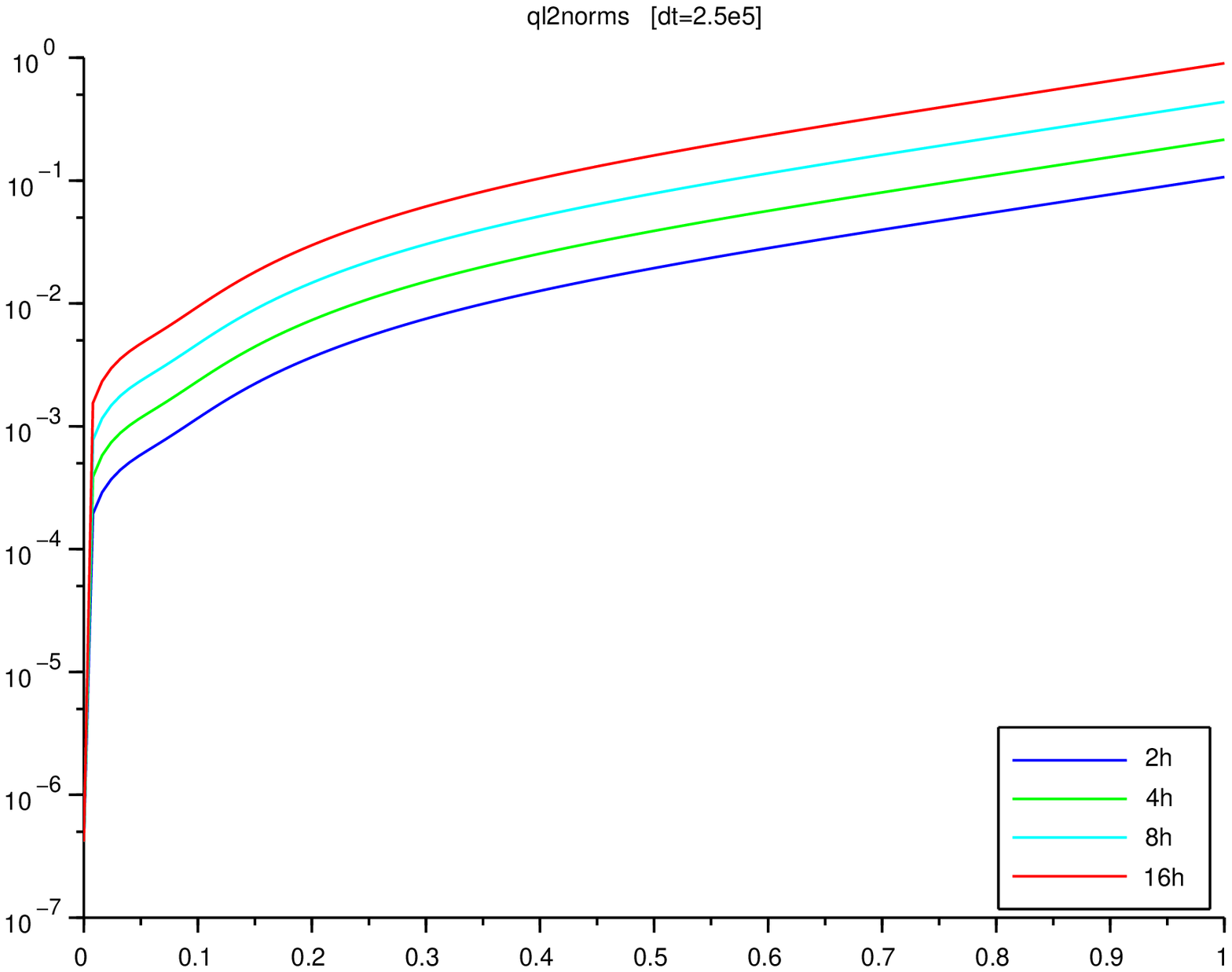}
\end{subfigure}%
\begin{subfigure}{.5\textwidth}
\centering
\includegraphics[width=70mm]{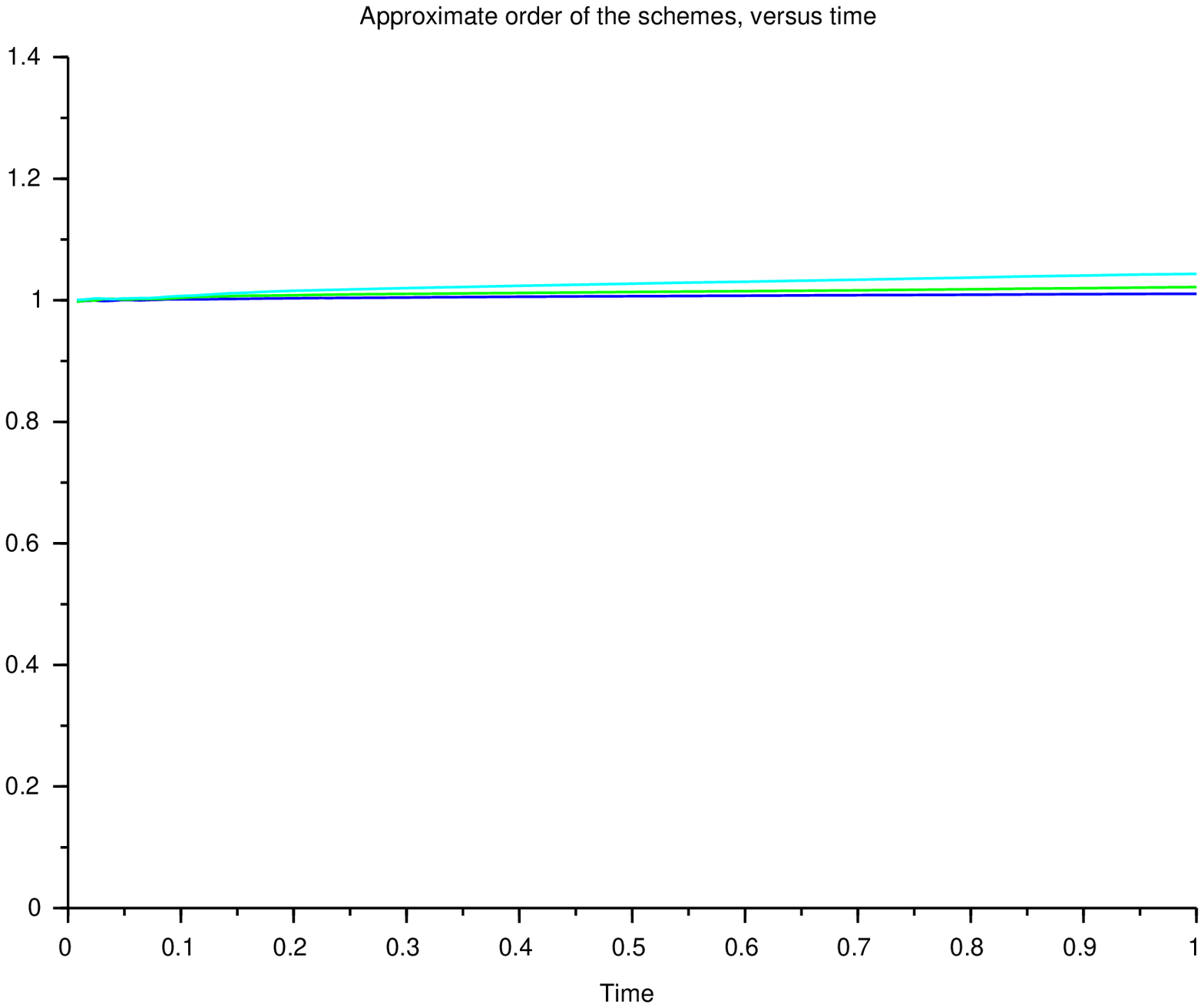}
\end{subfigure}%
\caption{\label{figEE}For the Euler-Maruyama method applied to NELD, we show the
$\ell^2$ norm of the position differences (left) and estimated order of convergence graph as a function of time (right), which exhibits first-order convergence, consistent with the error analysis.}
\end{figure}

The Euler-Maruyama integrator for NELD differs from the SE-A algorithm above
since we do not need to update the position before integrating the momentum,
leading to only one application of the periodic boundary conditions.  Thus we
get the algorithm:
\begin{algorithm}[H]\small
\caption{Euler-Maruyama\label{algEE}}
\begin{algorithmic}
\Statex 
\For{$k \gets 1 \dots $ Nsteps}      
\State \Call{GenKR}{$\qb^k,\pb^k$}
\Let{$ q_{tmp} $}{$  \pb^k + \pb^k \Delta t $}
\Let{$ \pb^{k+1} $}{$  \pb^k +  \Call{F}{\qb^k,\pb^k} \Delta t+\sigma \sqrt{\Delta t} \etab$}
\Let{$ \qb^{k+1} $}{$ q_{tmp} $}
\EndFor
\end{algorithmic}  
\end{algorithm}
Writing the update rule, and taking into account the application of PBCs, we have
\begin{gather}\label{expEE}
\vspace{0.8 cm}
\begin{split}
\begin{bmatrix} \qb^{k+1} \\ \pb^{k+1} \end{bmatrix}
&=
\begin{bmatrix} \qb^{k}- L \nb \\ \pb^{k}-A L \nb \end{bmatrix}+
\begin{bmatrix}
\pb^k -A L \nb \\
\Call{F}{\pb^k -A L \nb,\qb^k - L \nb}  
\end{bmatrix} \Delta t 
+
\begin{bmatrix}
0&0\\
0&\sigma I
\end{bmatrix} \sqrt{\Delta t} \etab\\
&=
\begin{bmatrix} \qb^{k}- L \nb \\ \pb^{k}-A L \nb \end{bmatrix}+
\begin{bmatrix}
\pb^k -A L \nb \\
\Call{F}{\qb^k,\pb^k}-A^2 L \nb  
\end{bmatrix}\Delta t
+
\begin{bmatrix}
0&0\\
0&\sigma I
\end{bmatrix} \sqrt{\Delta t} \etab\\
\end{split}
\end{gather}
Comparing with~\eqref{eq:ito-tay}, we compute the leading orders of the 
local truncation error, finding $O(\Delta t^{3/2})$ stochastic terms
and $O(\Delta t^2)$ deterministic terms, 
\begin{gather}\label{q2}
T_2= 
\begin{bmatrix}
\Call{F}{\pb,\qb}\\
(- \Call{$\nabla^2 E$}{\qb} +\gamma A) \pb +(A-\gamma I)\Call{F}{\pb,\qb} 
\end{bmatrix} \frac{\Delta t^2}{2} 
+\sigma \Delta t^{3/2} \begin{bmatrix}
0 & I\\
0 &  (A-\gamma I)\\
\end{bmatrix} \left(\frac{1}{2}\etab+\frac{1}{2\sqrt{3}}\zetab\right) 
\end{gather}
Therefore, the scheme will converge to first order, which is confirmed 
numerically in Figure~\ref{figEE}.

\subsection{Symplectic Euler B (SE-B)}
In Symplectic Euler B, the momentum is integrated first, then the position.  
The periodic boundary conditions need only be applied a single time during the
inner loop, and we have the following pseudocode:
\begin{algorithm}[H]\small
\caption{Symplectic Euler B
\label{algSEB}}
\begin{algorithmic}
\Statex
\For{$k \gets 1 \dots $ Nsteps}
\State \Call{GenKR}{$p,q$} 
\Let{$ \pb^{k+1} $}{ $\pb^k$  + \Call{F}{$\pb^{k+1},\qb^k $} $dt+\sigma dW$}
\Let{$ \qb^{k+1} $}{$  \qb^k + \pb^{k+1} \Delta t$}
\EndFor
\end{algorithmic}  
\end{algorithm}
\begin{figure}
\includegraphics[width=70mm]{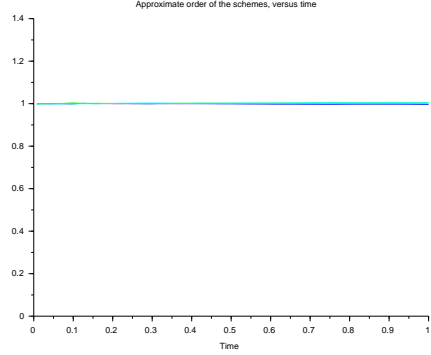}
\caption{\label{Opb0}Order of convergence in position, using SE-B}
\end{figure}
The numerical scheme is implicit in $\pb,$ though it is linear in $\pb^{k+1},$
and a perturbation of the identity, so that we can solve for $\pb^{k+1}$ and
expand in powers of $\Delta t,$ getting
\begin{equation*}
\begin{split}
\begin{bmatrix} \qb^{k+1} \\ \pb^{k+1} \end{bmatrix}
&=
\begin{bmatrix} \qb^{k}-L \nb \\ \pb^{k}-AL \nb \end{bmatrix}+
\begin{bmatrix}
\pb^k -AL \nb +(\Call{F}{\pb^{k+1},\qb^k - L \nb})\Delta t  \\
\Call{F}{\pb^{k+1},\qb^k - L \nb}  
\end{bmatrix}\Delta t+
\sigma\begin{bmatrix}
0 &0\\
0 &I
\end{bmatrix}dW\\
&=
(I+(\gamma I -A)\Delta t)^{-1} \left( \begin{bmatrix} 
(I+(\gamma I -A)\Delta t) (\qb^{k}- L \nb)  \\ \pb^{k}-A L \nb \end{bmatrix}  \right. \\
&\qquad + 
\begin{bmatrix}
\pb^k -A L \nb + (  - \nabla E(\qb^k)  +\gamma A(\qb^k- L \nb))  \Delta t  \\
(  - \nabla E(\qb^k)  +\gamma A(\qb^k- L \nb)) 
\end{bmatrix}\Delta t+
\sigma\begin{bmatrix}
0 & \Delta t\\
0 &I
\end{bmatrix} \sqrt{\Delta t} \etab \big) \\
&=
\begin{bmatrix} \qb^{k}- L \nb\\ \pb^{k}-A L \nb \end{bmatrix}+
\begin{bmatrix} \pb^k -A L \nb \\  \Call{F}{\pb^{k},\qb^k}-A^2 L \nb \end{bmatrix} \Delta t+
\begin{bmatrix}
 \Call{F}{\pb^{k},\qb^k}-A^2 L \nb  \\
(A-\gamma I)( \Call{F}{\pb^{k},\qb^k}-A^2 L \nb )   
\end{bmatrix} \Delta t^2 \\
&\quad  +
\sigma  \begin{bmatrix}
0 &0\\
0 &I
\end{bmatrix}\Delta t^{1/2} \etab  +
\sigma  \begin{bmatrix}
0 &I\\
0 &(A-\gamma I)
\end{bmatrix}\Delta t^{3/2} \etab 
+ O(\Delta t^{5/2})\\
\end{split}
\end{equation*}
As in the case of the Euler-Maruyama scheme, we find that the local truncation error is $O(\Delta t^{3/2})$ in the 
stochastic terms and $O(\Delta t^2)$ in the deterministic terms, 
\begin{equation}\label{T_3}
T_3= 
\begin{bmatrix}
-\Call{F}{\pb,\qb}\\
(- \Call{$\nabla^2 E$}{\qb} +\gamma A) \pb -(A-\gamma I)\Call{F}{\pb,\qb} 
\end{bmatrix} \frac{\Delta t^2}{2} 
 +\sigma \Delta t^{3/2} \begin{bmatrix}
0 & I\\
0 &  (A-\gamma I)\\
\end{bmatrix} \left(-\frac{1}{2}\etab+\frac{1}{2\sqrt{3}}\zetab\right) 
\end{equation}
The global truncation error converges to the first order. The numerical results illustrated in Figure~\ref{Opb0}  confirm the analytical result.

\subsection{Symplectic  Euler A Corrected (SE-AC) and ABAPO Corrected (ABAPO-C)}
The difference between ABA-O/SE-AC and the corrected schemes ABAPO-C/SE-AC 
presented here resides in the fact that applying periodic boundary conditions is 
only done once
during the scheme, while interparticle forces  $-\Call{$\nabla E$}{\qb}$ are 
computed using the periodic conditions (while particles may rest outside of the box).
Thus we get algorithms~\ref{algSEAC} and~\ref{ABAPO-C} for SE-AC and ABAPO-C, respectively.
\begin{algorithm}[H]\small
\caption{Symplectic Euler A Corrected (SE-AC)
\label{algSEAC}}
\begin{algorithmic}
\Statex
\For{$k \gets 1 \dots $ Nsteps}
\State \Call{GenKR}{$\qb^k,\pb^k$}  
\Let{$ \qb^{k+1} $}{$  \qb^k + \pb^k \Delta t$}  
\Let{$ \pb^{k+1} $}{$  \pb^k +  \Call{${F}_{\rm{PBC}}$}{$ $\pb^k,$ $\qb^{k+1} $ $} \Delta t + \sigma \sqrt{\Delta t} \etab$}
\EndFor
\end{algorithmic}  
\end{algorithm}
\begin{algorithm}[H]\small
\caption{ABAPO-C
\label{ABAPO-C}} 
\begin{algorithmic}
\Statex

\For{$k \gets 1 \dots $ Nsteps}
\State \Call{GenKR}{$\qb^k,\pb^k$} \Comment{Apply PBCs at the beginning of each iteration }
\Let{$ \pb^{k+\frac{1}{4}} $}{$  \pb^k + \frac{\Delta t}{2} \Call{$\nabla E$}{\qb^k}  $}
\Let{$ \qb^{k+1} $}{$  \qb^k +  \frac{\Delta t}{2} \pb^{k+\frac{1}{4}} $}
\Let{$  \pb^{k+ \frac{1}{2}} $}{$\pb^{k+ \frac{1}{4}} + \frac{\Delta t}{2} {\Call{$\nabla E_{\rm{PBC}}$}{ \qb^{k+1}}} $}
\Let{$ \tilde{ \pb^{k+ \frac{1}{2}}}   $}{$ \exp(\Delta t A) \pb^{k+\frac{1}{2}} $}
\Let{$ \pb^{k+1} $}{$ \gamma \tilde{ \pb^{k+ \frac{1}{2}}}  +(1-\gamma)A \qb^{k+1} +\sigma dW   $}
\EndFor
\end{algorithmic}  
\end{algorithm}

\begin{figure}
\includegraphics[width=70mm]{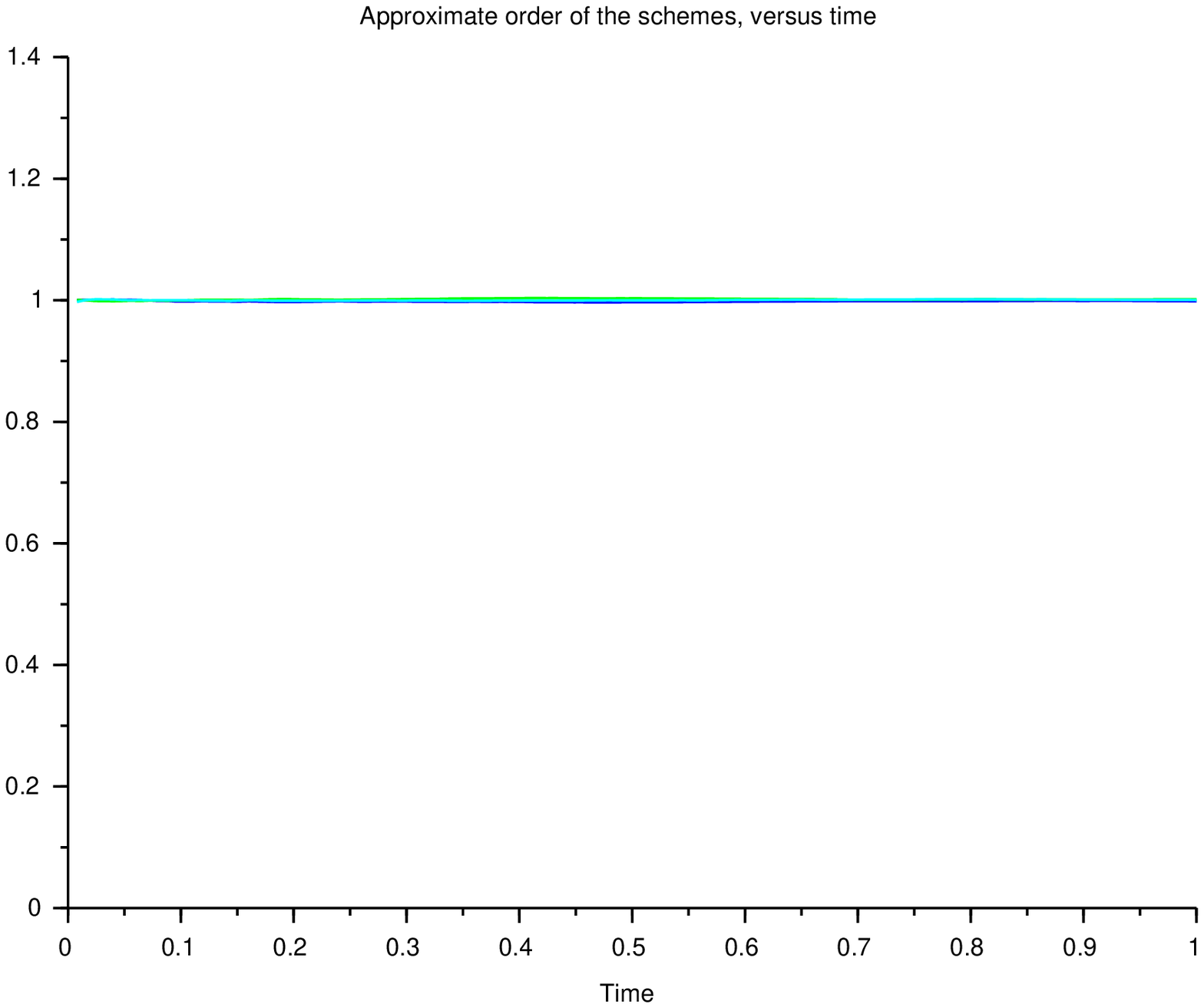}
\includegraphics[width=70mm]{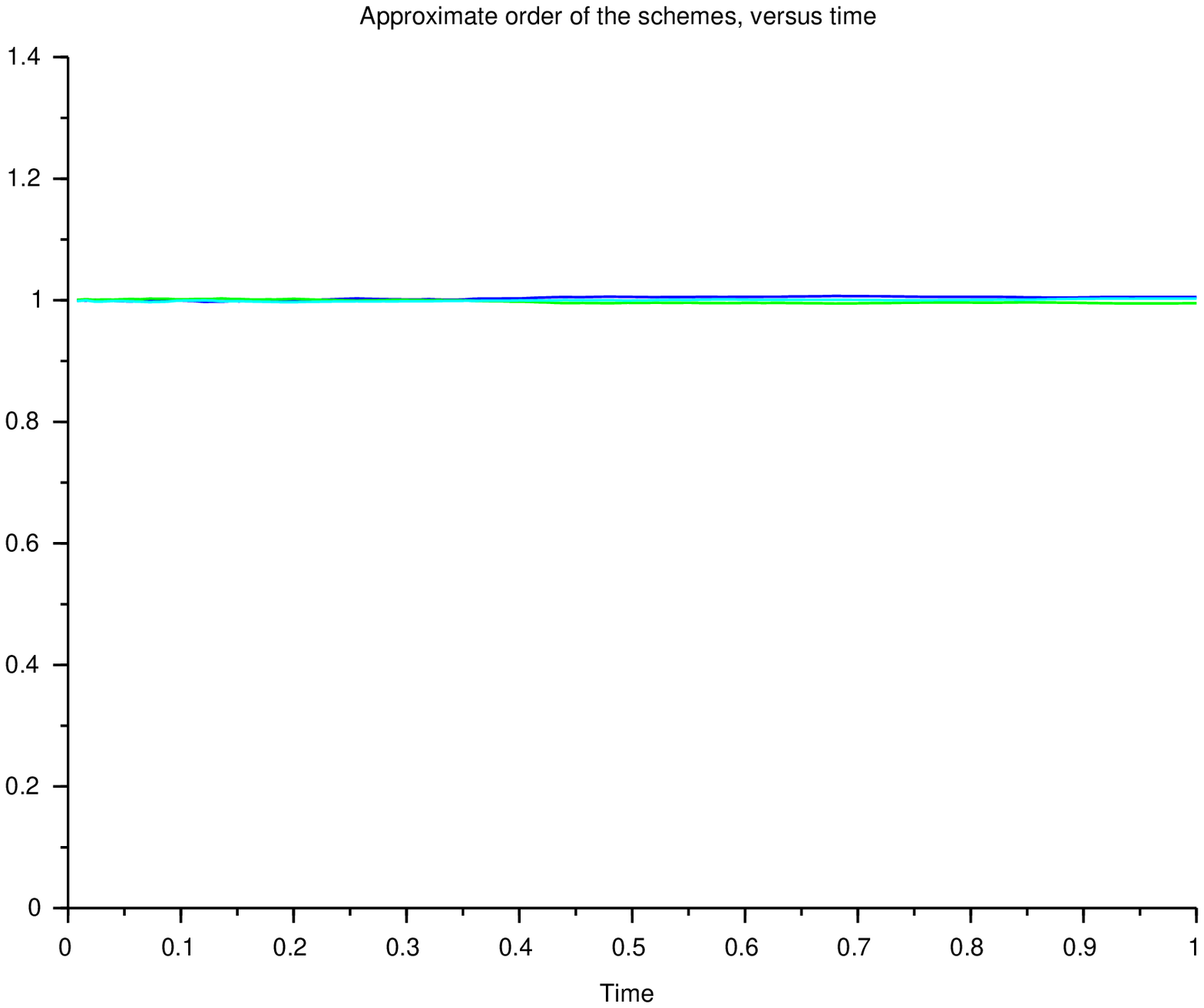}
\caption{\label{figABAPO-C}Order of convergence in position, using SE-AC (left) and ABAPO-C (right)}
\end{figure}
Expanding out SE-AC, we have a similar expression to SE-A~\eqref{eq:sea_exp}, though there
is only a single application of periodic boundary conditions, so we have
\begin{equation}
\begin{split}
\begin{bmatrix} \qb^{k+1} \\ \pb^{k+1} \end{bmatrix}
&=
\begin{bmatrix} \qb^{k}-L \nb \\ \pb^{k}-AL\nb \end{bmatrix}+
\begin{bmatrix}
\pb^k -AL\nb \\
\Call{F}{\pb^k -AL\nb,\qb^k -L\nb+(\pb^k -AL\nb ) \Delta t}  
\end{bmatrix} \Delta t+
\begin{bmatrix}
0&0\\
0&\sigma I
\end{bmatrix} \etab\\
&=
\begin{bmatrix} \qb^{k}-L\nb \\ \pb^{k}-AL\nb \end{bmatrix} 
+
\begin{bmatrix}
\pb^k -AL\nb \\
\Call{F}{\qb^k,\pb^k} -A^2L \nb  
\end{bmatrix} \Delta t
+
\begin{bmatrix}
0&0\\
0&\sigma I
\end{bmatrix} \sqrt{\Delta t} \etab\\
&
\quad +
\begin{bmatrix}
0 \\
\gamma A(\pb^k-AL\nb) - \Call{$\nabla^2 E$}{\qb^k}(\pb^k -AL \nb )
\end{bmatrix} \Delta t^2.
\end{split}
\end{equation}
Then the leading order terms in the local truncation error are $O(\Delta t^{3/2})$ in
the stochastic terms and $O(\Delta t^2)$ in the deterministic terms,
\begin{equation*}\label{T4}
\begin{split}
T_4 &=
\begin{bmatrix}
\Call{F}{\pb,\qb} + A^2 L \nb\\
(\Call{$\nabla^2 E$}{\qb} -\gamma A) \, \pb +(A-\gamma I)\Call{F}{\pb,\qb} 
\end{bmatrix} \frac{\Delta t^2}{2} 
+\sigma \Delta t^{3/2} \begin{bmatrix}
0 & I\\
0 &  (A-\gamma I)\\
\end{bmatrix} \left(\frac{1}{2}\etab+\frac{1}{2\sqrt{3}}\zetab\right) 
\end{split}
\end{equation*}
Therefore, we expect to find first order convergence for the scheme, which
is precisely what we observe in Figure~\ref{figABAPO-C}.  We similarly see the 
same improvement in convergence for the ABAPO-C scheme.

\section{Second Order Integrator of the Langevin Equation A and B (SOILE-A \& B)}
We base our Second order NELD integrators on algorithms developed for equilibrium Langevin 
dynamics in~\cite{Eric}. Since we need to integrate the position first in both methods, we 
apply the ideas from the corrected algorithm where we wait to remap the particle positions
until the end of the variable updates.  The standard SOILE-A scheme is a generalization of 
the Langevin equation for velocity-Verlet algorithm while, SOILE-B is a quasi-symplectic
scheme.  
The algorithms are described as follows:
\begin{algorithm}
\caption{Second Order Integrator of the Langevin Equation A (SOILE-A)
\label{algSOILE-A}}
\begin{algorithmic}
\Statex
\For{$k \gets 1 \dots $ Nsteps}
\State \Call{GenKR}{$\qb^k,\pb^k$}  
\Let{$ \qb^{k+1} $}{$  \qb^k + \pb^k \Delta t+ \Call{$F$}{$ $\pb^k,$ $\qb^{k} $ $} \frac{\Delta t^2}{2} + \Delta t^{3/2}(\frac{1}{2}\etab^k+\frac{1}{2\sqrt{3}}\zetab^k) $}  
\Let{$ \pb^{k+1} $}{$  \pb^k +  (\Call{${F}_{\rm{PBC}}$}{$ $\pb^k,$ $\qb^{k+1} $ $}+ \Call{$F$}{$ $\pb^k,$ $\qb^{k} $ $})\frac{\Delta t^2}{2}+\sigma dW-(\gamma I - A)(\Call{$F$}{$ $\pb^k,$ $\qb^{k} $ $} \frac{\Delta t^2}{2} + \Delta t^{3/2}(\frac{1}{2}\etab^k+\frac{1}{2\sqrt{3}}\zetab^k) )$}
\EndFor
\end{algorithmic}
\end{algorithm}

\begin{figure}[h]
\begin{subfigure}{.5\textwidth}
\centering
\includegraphics[width=70mm]{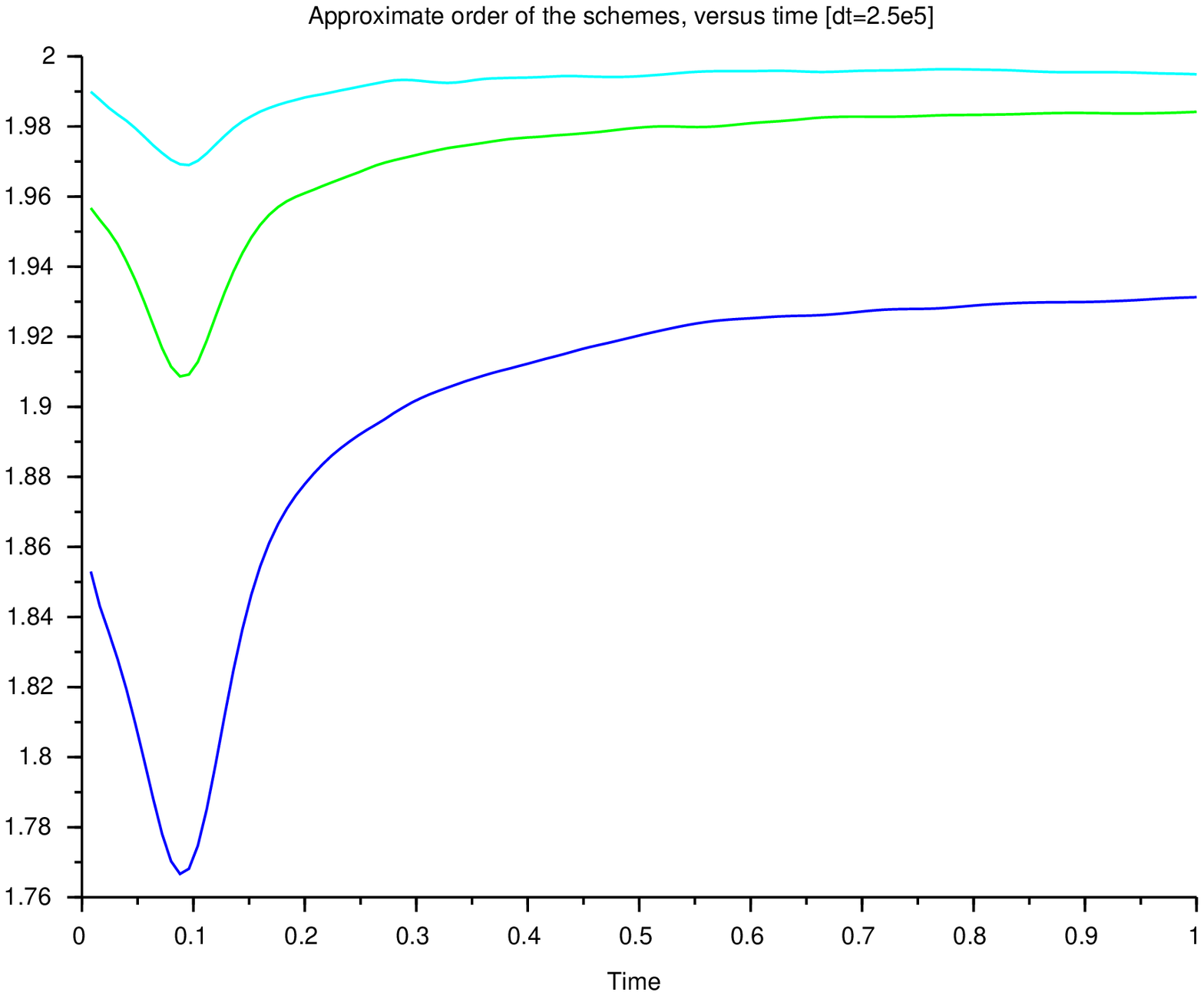}
\end{subfigure}%
\begin{subfigure}{.5\textwidth}
\centering
\includegraphics[width=70mm]{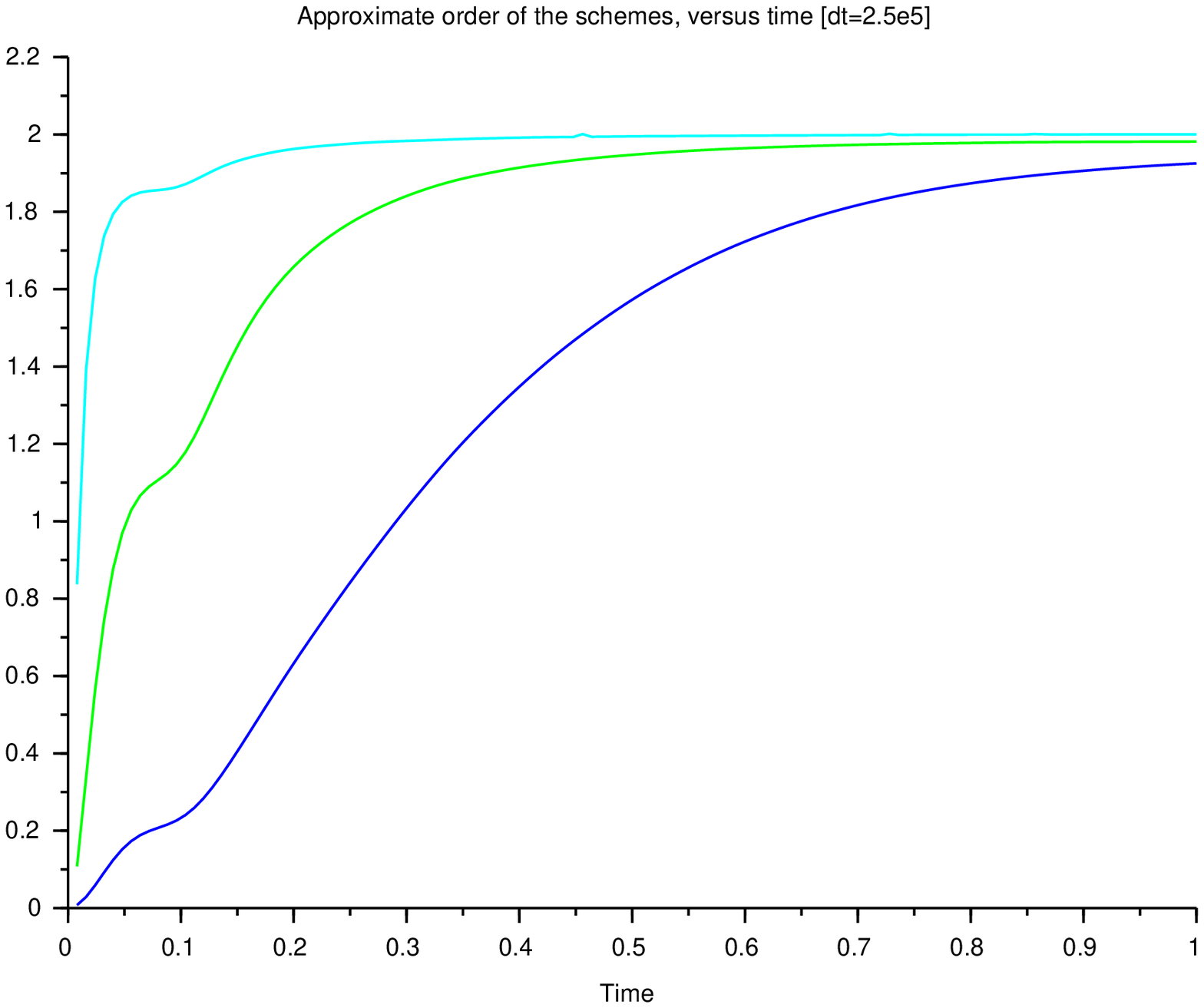}
\end{subfigure}%
\caption{\label{figSOILE}Estimated error of convergence for the SOILE-A (left) and SOILE-B (right) schemes.
Both shown to be second order based on truncation error analysis, and the observed errors
order is approximately second order.  Note the different axes.}
\end{figure}

\begin{algorithm}
\caption{Second Order Integrator of the Langevin Equation B (SOILE-B)
\label{SEA}}
\begin{algorithmic}
\Statex
\For{$k \gets 1 \dots $ Nsteps}
\State \Call{GenKR}{$\qb^k,\pb^k$} 
\Let{$ \pb^{k+1/2} $}{$  \pb^k +  \frac{1}{2}\Call{$F$}{$ $\pb^k,$ $\qb^{k} $ $}\Delta t 
-\frac{1}{4}(\gamma I - A) \big[ \Call{$F$}{$ $\pb^k,$ $\qb^{k} $ $}\frac{1}{2}\Delta t^2
+\sigma\Delta t^{\frac{3}{2} }(\frac{1}{2}\etab +\frac{1}{\sqrt{3} }\zetab)} \big] $
\Let{$ \qb^{k+1} $}{$  \qb^k + \pb^{k+1/2} \Delta t+ \sigma \Delta t^{\frac{3}{2} }\frac{1}{\sqrt{3} }\zetab$}
\Let{$ \pb^{k+1} $}{$  \pb^{k+1/2} +  \frac{1}{2} \Call{$\tilde{F}_{\rm{PBC}}$}{ $
$\pb^{k+1/2},$ $\qb^{k+1} $ $}\Delta t + \frac{1}{2} \sigma \sqrt{\Delta t} \etab
-\frac{1}{4}(\gamma I - A) \big[ \Call{$\tilde{F}_{\rm{PBC}}$}{$ $\pb^{k+1/2},$
$\qb^{k+1} $ $}\frac{1}{2}\Delta t^2+\sigma\Delta t^{\frac{3}{2}
}(\frac{1}{2}\etab +\frac{1}{\sqrt{3} }\zetab)} \big] $ 
\EndFor
\end{algorithmic}
\end{algorithm}
We write out the SOILE-A scheme, and expand out to second order, giving 
\begin{gather}\label{extSOILE-A}
\begin{split}
\tilde{ \begin{bmatrix} \qb^{k+1} \\ \pb^{k+1} \end{bmatrix}}
&=
\begin{bmatrix} \qb^{k}- L \nb \\ \pb^{k}-A L \nb \end{bmatrix}+
\begin{bmatrix}
\pb^k -A L \nb \\
  \frac{1}{2}(-\Call{$\nabla E$}{\qb^{k+1}}-\Call{$\nabla E$}{\qb^k}  + \gamma A(\qb^{k+1}+\qb^{k}))-(\gamma I - A)(\pb^k-A L \nb)
\end{bmatrix}\Delta t\\&\quad
+\begin{bmatrix}
0&0\\
0&\sigma I
\end{bmatrix} \sqrt{\Delta} \etab
+\begin{bmatrix}
\Call{F}{\qb^k,\pb^k}-A^2 L \nb \\
0
\end{bmatrix} \frac{\Delta t^2}{2} 
 +
\begin{bmatrix}
0 & \sigma I\\
0 & \sigma (A-\gamma I)\\
\end{bmatrix} \Delta t^{3/2}(\frac{1}{2}\zetab^k+\frac{1}{2\sqrt{3}}\etab^k)\\
&=
\begin{bmatrix} \qb^{k}- L \nb \\ \pb^{k}-A L \nb \end{bmatrix}+
\begin{bmatrix}
\pb^k -A L \nb \\
\Call{F}{\qb^k,\pb^k}-A^2 L \nb 
\end{bmatrix}\Delta t+
\begin{bmatrix}
0&0\\
0&\sigma I
\end{bmatrix} \sqrt{\Delta} \etab\\& \quad
+\begin{bmatrix}
\Call{F}{\qb^k,\pb^k}-A^2 L \nb \\
-( \Call{$\nabla^2 E$}{\qb^k} +\gamma A)(\pb^k -A L \nb) +(A-\gamma I)(\Call{F}{\qb^k,\pb^k}-A^2 L \nb ) 
\end{bmatrix} \frac{\Delta t^2}{2} \\
&\quad +
\begin{bmatrix}
0 & \sigma I\\
0 & \sigma (A-\gamma I)\\
\end{bmatrix} \Delta t^{3/2} \left(\frac{1}{2}\zetab^k+\frac{1}{2\sqrt{3}}\etab^k\right)
+O(\Delta t^3)
\end{split}
\end{gather}
These terms cancel with the exact expansion~\eqref{eq:ito-tay}, with $O(\Delta t^{5/2})$ 
stochastic terms and $O(\Delta t^3)$ deterministic terms, giving a globally second-order
convergent scheme.  We plot the order of convergence for both schemes in Figure~\ref{figSOILE},
where we observe    
Moreover, we observe that the graphs of convergence~\ref{figSOILE} satisfy the second order scheme as predicted the truncation error.  Note that both exhibit lower initial convergence, leveling off close to second order.  

\section{Conclusion}

We have derived several numerical integrators for nonequilibrium Langevin dynamics and have shown that care
must be taken in applying the periodic boundary conditions, or there can be a breakdown in the order of
convergence.  Provided that the pbcs are not applied in the middle of an update step, we have demonstrated
several prototypical schems of order one and two applied to NELD.  For these orders, deforming the domain is 
performed after all other updates, and we still acheive the desired accuracy.

Several extensions are possible.  First, deriving conditions for general higher-order schemes that 
appropriately incorporate the deforming simulation box and nonequilibrium PBCs, for example, for 
general stochastic Runge-Kutta schemes or variational schemes will be of interest.  Also, of large
interest in molecular dynamics is the convergence to the invariant measure as in~\cite{Burrage, Leimkuhler}.  
This is challenging in the present case, since unlike Langevin dynamics, there is not generally
an analytic expression for the invariant measure of the original dynamics~\eqref{neld}.  

\section{Acknowledgements}

MD was supported by the DARPA EQUiPS program.

\end{document}